\documentclass[a4paper, 12pt]{article}

\usepackage[english]{babel}

\usepackage[top=2cm,bottom=2cm,left=3cm,right=3cm,marginparwidth=1.75cm]{geometry}

\usepackage{amsmath}
\usepackage{graphicx}
\usepackage[colorlinks=true, allcolors=blue]{hyperref}
\newtheorem{theorem}{Theorem}
\newtheorem{lemma}{Lemma}
\newtheorem{definition}{Definition}
\newtheorem{corollary}{Corollary}

\title{Minimax approach to the estimation problem for homogeneous random fields}
\author{Oleksandr Masyutka, Mikhail Moklyachuk \\ Taras Shevchenko National University of Kyiv \\ omasyutka@gmail.com, moklyachuk@gmail.com}
\date{}

\begin{document}
\maketitle

\begin{abstract}
The problem  of the mean-square optimal estimation of the linear functionals  which depend on the unknown values of a multidimensional homogeneous random field from observations of the field with noise is considered. The minimax (robust) method of estimation is applied in the case where the spectral densities of the fields are not known exactly while some sets of admissible spectral densities are given. Formulas that determine the least favourable spectral densities and the minimax spectral characteristics are derived for some special sets of admissible densities.
\end{abstract}

\section{Introduction}

Traditional methods of solution of the linear extrapolation, interpolation and filtering problems for
stationary stochastic processes and homogeneous random fields are developed
 under the condition that spectral densities  of processes and fields are known exactly
 (see, for example, selected works of A.~N.~Kolmogorov \cite{Kolmogorov}, 
 survey by T.~Kailath \cite{Kailath}, 
 books by Yu.~A.~Rozanov \cite{Rozanov}, E.~J.~ Hannan \cite{Hannan},
  A.~Malyarenko \cite{Malyarenko},
V.~S.~Mandrekar and D.~A.~Redett \cite{Mandrekar}, N.~Wiener \cite{Wiener}, M.~I.~Yadrenko \cite{Yadrenko},
 A.~M.~Yaglom \cite{Yaglom1961,Yaglom1987}).

 The basic assumption of most of the methods of  estimation of the unobserved values of stochastic processes and random fields is that the spectral densities of the considered stochastic processes and random fields are exactly known.
However, in practice, these methods are not applicable since the complete information on the spectral densities is impossible in most cases.
In order to solve the problem  parametric or nonparametric estimates of the unknown spectral densities are found. 
Then, one of  traditional estimation methods is applied, provided that the selected densities are the true ones. This procedure can result in significant increasing of the value of error as K.~S.~Vastola and H.~V.~Poor \cite{Vastola} have demonstrated with the help of some examples.
To avoid this effect one can search the estimates which are optimal for all densities from a certain class of admissible spectral densities. 
These estimates are called minimax since they minimize the maximum value of the error. 
The paper by Ulf~Grenander \cite{Grenander} should be marked as the first one where this approach to extrapolation problem for stationary processes was proposed.
Several models of spectral uncertainty and minimax-robust methods of data processing can be found in the survey paper by S.~A.~ Kassam and H.~V.~ Poor \cite{Kassam}. 
In the paper by J.~Franke \cite{Franke} the minimax extrapolation problem for stationary sequences is investigated with the help of convex optimization methods.
This approach makes it possible to find equations that determine the least favorable spectral densities for various classes of densities.
 In the book by M.~ Moklyachuk and A.~ Masyutka  a  minimax  technique of the estimation for vector-valued stationary stochastic  processes is proposed \cite{Masyutka}. 
 Estimation problems for random fields was considered in the paper by M. Moklyachuk and N. Shchestyuk \cite{Shestyuk}.
In the book by M.~P.~Moklyachuk, O.~Yu.~Masyutka, I.~I.~Golichenko \cite{Moklyachuk2018}
the minimax approach is applied to investigate the estimation problems
for functionals which depend on the unknown values of stationary random fields on a sphere.

In this paper we deal with the problem of the mean-square optimal linear estimation of the functionals
\[A_K\vec{\xi}=\iint\limits_{K}\vec{a}(s,t)^\top\vec{\xi}(s,t)dsdt,\]
which depend on the unknown values of a multidimensional homogeneous random field $\vec{\xi}(s,t)=\{\xi_k(s,t)\}_{k=1}^{N}$ from observations of the field $\vec{\xi}(s,t)+\vec{\eta}(s,t)$ at points  $(s,t)\in R^2\setminus K$, where $\vec{\eta}(s,t)=\{\eta_k(s,t)\}_{k=1}^{N}$ is an uncorrelated with $\vec{\xi}(s,t)$ multidimensional homogeneous random field. The case  of spectral uncertainty is considered. 
Formulas for calculating the spectral characteristic $h(F,G)$ and the mean square error $\Delta(F,G)$ of the optimal linear estimate of the functionals under the condition that spectral densities $F(\lambda,\mu), G(\lambda,\mu)$ of the fields are exactly known were derived in  \cite{Moklyachuk2026}.
In the case of spectral uncertainty, where the spectral densities are not exactly known while a set of admissible spectral densities is given, the minimax method is applied. 
Formulas for determination the least favorable spectral densities and the minimax-robust spectral characteristics of the optimal estimates of the functionals are proposed for some specific classes of admissible spectral densities.

\section{Minimax method of interpolation}

The formulas proposed in \cite{Moklyachuk2026} for calculating the spectral characteristic $h(F,G)$ and the mean square error $\Delta(F,G)$ of the optimal linear estimate of the functionals may be employed under the condition that spectral densities $F(\lambda,\mu), G(\lambda,\mu)$ of the fields are exactly known. In the case where the densities are not known exactly while a set $D = D_{F} \times D_{G} $ of possible spectral densities is given, the minimax (robust) approach to estimation of functionals of the unknown values of multidimensional homogeneous random fields is reasonable. 
Instead of searching an estimate that is optimal for a given spectral densities  we find an estimate that minimizes the mean square error for all spectral densities $F(\lambda,\mu), G(\lambda,\mu)$ from a given class $D_F\times D_G$ simultaneously.

 \begin{definition}
 For a given class of spectral densities $D=D_F \times D_G$ the spectral densities  $F^0(\lambda,\mu) \in D_F$, $G^0(\lambda,\mu) \in D_G$ are called least favorable in the class $D$ for the optimal linear estimation of the functional $A_K\vec{\xi}$  if the following relation holds true $$\Delta\left(F^0,G^0\right)=\Delta\left(h\left(F^0,G^0\right);F^0,G^0\right)=\max\limits_{(F,G)\in D_F\times D_G}\Delta\left(h\left(F,G\right);F,G\right).$$
\end{definition}

\begin{definition}
For a given class of spectral densities $D=D_F \times D_G$ the spectral characteristic $h^0(\lambda,\mu)$ of the optimal linear estimation of the functional $A_K\vec{\xi}$ is called minimax-robust if there are satisfied conditions
$$h^0(\lambda,\mu)\in H_D= \bigcap\limits_{(F,G)\in D_F\times D_G} L_2^{K-}(F+G),$$
$$\min\limits_{h\in H_D}\max\limits_{(F,G)\in D}\Delta\left(h;F,G\right)=\max\limits_{(F,G)\in D}\Delta\left(h^0;F,G\right),$$
where $L_2^{K-}(F+G)$ is the subspace of the space $L_2(F+G)$ defined in \cite{Moklyachuk2026}.
\end{definition}

The least favorable spectral densities $F^0(\lambda,\mu)$, $G^0(\lambda,\mu)$ and the minimax spectral characteristic $h^0=h(F^0,G^0)$ form a saddle point of the function $\Delta \left(h;F,G\right)$ on the set  $H_D\times D.$ The saddle point inequalities
\[\Delta\left(h^0;F,G\right)  \leq\Delta\left(h^0;F^0,G^0\right)\leq \Delta\left(h;F^0,G^0\right), \quad \forall h \in H_D, \forall F \in D_F, \forall G \in D_G,\]
hold true if  $h^0=h(F^0,G^0)$ and $h(F^0,G^0)\in H_D,$ where $(F^0,G^0)$ is a solution to the constrained optimization problem
\begin{equation} \label{7}
\sup\limits_{(F,G)\in D_F\times D_G}\Delta\left(h(F^0,G^0);F,G\right)=\Delta\left(h(F^0,G^0);F^0,G^0\right),
\end{equation}
\[\Delta \left(h\left(F^{0} , G^{0} \right); F, G\right)=\]
\[=\frac{1}{4\pi^2}\int_{-\infty}^{\infty}\int_{-\infty}^{\infty}((A_{K} (\lambda,\mu ))^\top G^0(\lambda,\mu )+(C_{K}^0(\lambda,\mu ))^\top)(F^0(\lambda,\mu )+G^0(\lambda,\mu ))^{-1} F(\lambda,\mu )\times\]
\[\times (F^0(\lambda,\mu )+G^0(\lambda,\mu ))^{-1}((A_{K} (\lambda,\mu ))^\top G^0(\lambda,\mu )+(C_{K}^0(\lambda,\mu ))^\top)^*d\lambda d\mu+\]
\[+\frac{1}{4\pi^2}\int_{-\infty}^{\infty}\int_{-\infty}^{\infty}((A_{K} (\lambda,\mu ))^\top F^0(\lambda,\mu )-(C_{K}^0(\lambda,\mu ))^\top)(F^0(\lambda,\mu )+G^0(\lambda,\mu ))^{-1} G(\lambda,\mu )\times\]
\[\times (F^0(\lambda,\mu )+G^0(\lambda,\mu ))^{-1}((A_{K} (\lambda,\mu ))^\top F^0(\lambda,\mu )-(C_{K}^0(\lambda,\mu ))^\top)^*d\lambda d\mu,\]
where $A_{K} (\lambda,\mu )$, $C_{K} (\lambda,\mu )$ were defined in \cite{Moklyachuk2026}.

From the introduced definitions and formulas derived in  \cite{Moklyachuk2026} we can obtain the following statements.

\begin{lemma}
Spectral densities $F^0(\lambda,\mu)\in D_F$, $G^0(\lambda,\mu) \in D_G$, satisfying the minimality condition \cite{Moklyachuk2026} are the least favorable in the class  $D=D_F\times D_G$ for the optimal linear interpolation of the functional $A_K\vec{\xi}$ if the functions
\[(F^0(\lambda,\mu)+G^0(\lambda,\mu))^{-1}, F^0(\lambda,\mu)(F^0(\lambda,\mu)+G^0(\lambda,\mu))^{-1},\]
determine
\[C_K^0(\lambda,\mu)=\iint\limits_{K}\vec{c}(s,t)e^{i(s\lambda+t\mu)}dsdt,\]
which gives a solution to the constrained optimization problem \eqref{7}. The minimax spectral characteristic  $h^0=h(F^0,G^0)$ can be calculated by the formula \cite{Moklyachuk2026}
\[(h^0(\lambda,\mu))^\top=(A_K(\lambda,\mu))^\top F^0(\lambda,\mu)(F^0(\lambda,\mu)+G^0(\lambda,\mu))^{-1}-\]
\begin{equation} \label{sphar}
-(C_K^0(\lambda,\mu))^\top(F^0(\lambda,\mu)+G^0(\lambda,\mu))^{-1}.
\end{equation}
if $h(F^0,G^0) \in H_D$.
\end{lemma}

\begin{corollary}
Spectral densities $F^0(\lambda,\mu)\in D_F,$ $G^0(\lambda,\mu) \in D_G$ satisfying the minimality condition \cite{Moklyachuk2026} are the least favorable in the class  $D=D_F\times D_G$ for the optimal linear estimation of the functional $A_K\vec{\xi}$ where $K=R\times[0,T]$ if the functions
\[(F^0(\lambda,\mu)+G^0(\lambda,\mu))^{-1},\quad F^0(\lambda,\mu)(F^0(\lambda,\mu)+G^0(\lambda,\mu))^{-1},\]
\[F^0(\lambda,\mu)(F^0(\lambda,\mu)+G^0(\lambda,\mu))^{-1}G^0(\lambda,\mu)\]
determine the operators $\mathbf{B}_T^0, \mathbf{R}_T^0, \mathbf{Q}_T^0$ \cite{Moklyachuk2026}, which give a solution to the constrained optimization problem
\[
\max\limits_{(F,G)\in D_F\times D_G}\int_{-\infty}^{\infty}\langle(\mathbf{R}_T\mathbf{a})(\lambda,t),(\mathbf{B}_T^{-1}\mathbf{R}_T\mathbf{a})(\lambda,t)\rangle+\langle(\mathbf{Q}_T\mathbf{a})(\lambda,t),\vec{a}(\lambda,t)\rangle d\lambda=
\]
\begin{equation} \label{extrem}
=\int_{-\infty}^{\infty}\langle(\mathbf{R}_T^0\mathbf{a})(\lambda,t),((\mathbf{B}_T^0)^{-1}\mathbf{R}_T^0\mathbf{a})(\lambda,t)\rangle+\langle(\mathbf{Q}_T^0\mathbf{a})(\lambda,t),\vec{a}(\lambda,t)\rangle d\lambda.
\end{equation}
The minimax spectral characteristic  $h^0=h(F^0,G^0)$ is determined by the formula (\ref{sphar}) if $h(F^0,G^0) \in H_D.$
\end{corollary}

\begin{corollary}
Let the spectral density $F^0(\lambda,\mu)\in D_F$  satisfy the minimality condition \cite{Moklyachuk2026}.  The spectral density $F^0(\lambda,\mu)\in D_F$ is the least favorable in the class $ D_F$ for the optimal linear estimation of the functional $A_K\vec{\xi}$ from the observation of the field $\vec{\xi}(s,t)$ at points $(s,t)\in R^2\setminus(R\times[0,T])$ if the function $(F^0(\lambda,\mu))^{-1}$ determine the operator $\mathbf{B}_T^0$ \cite{Moklyachuk2026} which determines a solution to the constrain optimization problem
 \begin{equation} \label{extrem2}
\max\limits_{F\in D_F}\int_{-\infty}^{\infty}\langle (\mathbf{B}_T^{-1}\mathbf{a})(\lambda,t),\vec{a}(\lambda,t)\rangle d\lambda=\int_{-\infty}^{\infty}<((\mathbf{B}_T^0)^{-1}\mathbf{a})(\lambda,t),\vec{a}(\lambda,t)>d\lambda.
\end{equation}
The minimax spectral characteristic $h^0=h_T(F^0)$ is determined by the formula \cite{Moklyachuk2026} 
\begin{equation} \label{sphar1}
(h_T(F^0))^\top=(A_T(\lambda,\mu))^\top-(C_T^0(\lambda,\mu))^\top(F^0(\lambda,\mu))^{-1},
\end{equation}
if $h_T(F^0) \in H_{D_F}.$
 \end{corollary}

The constrained optimization problem (\ref{7}) is equivalent to the unconstrained optimization problem \cite{Pshenichnyj}:
\begin{equation} \label{8}
\Delta_D(F,G)=-\Delta(h(F^0,G^0);F,G)+\delta((F,G)\left|D_F\times D_G\right.)\rightarrow \inf,
\end{equation}
where $\delta((F,G)|D_F\times D_G)$ is the indicator function of the set $D=D_F\times D_G$. Solution of the problem (\ref{8}) is characterized by the condition $0 \in \partial\Delta_D(F^0,G^0),$ where $\partial\Delta_D(F^0,G^0)$ is the subdifferential of the convex functional $\Delta_D(F,G)$ at point $(F^0,G^0)$ \cite{Rockafellar}. This condition makes it possible to find the least favourable spectral densities in some special classes of spectral densities \cite{Ioffe}, \cite{Pshenichnyj}.

\begin{lemma}
Let $(F^0,G^0)$ be a solution to the optimization problem (\ref{8}). The spectral densities  $F^0(\lambda)$, $G^0(\lambda)$ are the least favorable in the class $D=D_F\times D_G$ and the spectral characteristic  $h^0=h(F^0,G^0)$ is the minimax of the optimal linear estimate of the functional $A_K\vec{\xi}$ if $h(F^0,G^0) \in H_D$.
\end{lemma}

\section{Least favorable spectral densities in the class \texorpdfstring{$D = D_{V}^U\times D_{2\delta}$}{}}

Consider the problem of  mean square optimal estimation of the functional $A_K\vec{\xi}$ where $K=R\times[0,T]$ in the case when spectral densities of the fields  belong to the class of admissible spectral densities  $D=D_{V}^U \times D_{2\delta}$,
\[{D_{V}^{U}} ^{1}  =\biggl\{F(\lambda,\mu )\bigg|{\mathrm{Tr}}\, V(\lambda,\mu)\le {\mathrm{Tr}}\, F(\lambda,\mu )\le {\mathrm{Tr}}\, U(\lambda,\mu ),\]
\[\frac{1}{2\pi } \int _{-\infty }^{\infty }{\mathrm{Tr}}\,  F(\lambda,\mu)d\mu =p(\lambda) \biggr\},\]
\[D_{2\delta}^{1}=\left\{G(\lambda,\mu )\biggl|\frac{1}{2\pi } \int_{-\infty }^{\infty }\left|{\rm{Tr}}(G(\lambda,\mu )-G_{1} (\lambda,\mu))\right|^{2} d\mu \le \delta(\lambda)\right\};\]
\[{D_{V}^{U}} ^{2}  =\biggl\{F(\lambda,\mu )\bigg|v_{kk} (\lambda,\mu )  \le f_{kk} (\lambda,\mu )\le u_{kk} (\lambda,\mu ),\]
\[\frac{1}{2\pi } \int _{-\infty }^{\infty }f_{kk} (\lambda,\mu)d\mu  =p_{k}(\lambda) , k=\overline{1,N}\biggr\},\]
\[D_{2\delta}^{2}=\left\{G(\lambda,\mu )\biggl|\frac{1}{2\pi } \int_{-\infty }^{\infty }\left|g_{kk} (\lambda,\mu )-g_{kk}^{1} (\lambda,\mu)\right|^{2} d\mu  \le \delta_{k}(\lambda), k=\overline{1,N}\right\};\]
\[{D_{V}^{U}} ^{3}  =\biggl\{F(\lambda,\mu )\bigg|\left\langle B_{1},V(\lambda,\mu )\right\rangle \le \left\langle B_{1} ,F(\lambda,\mu)\right\rangle \le \left\langle B_{1} ,U(\lambda,\mu )\right\rangle,\]
\[\frac{1}{2\pi }\int _{-\infty }^{\infty }\left\langle B_{1},F(\lambda,\mu)\right\rangle d\mu  =p(\lambda)\biggr\},\]
\[D_{2\delta}^{3}=\left\{G(\lambda,\mu )\biggl|\frac{1}{2\pi } \int_{-\infty }^{\infty }\left|\left\langle B_{2} ,G(\lambda,\mu )-G_{1}(\lambda,\mu )\right\rangle \right|^{2} d\mu  \le \delta(\lambda)\right\};\]
\[{D_{V}^{U}} ^{4}=\left\{F(\lambda,\mu )\bigg|V(\lambda,\mu )\le F(\lambda,\mu)\le U(\lambda,\mu ), \frac{1}{2\pi } \int _{-\infty }^{\infty}F(\lambda,\mu )d\mu=P(\lambda)\right\},\]
\[D_{2\delta}^{4}=\left\{G(\lambda,\mu )\biggl|\frac{1}{2\pi } \int_{-\infty }^{\infty }\left|g_{ij} (\lambda,\mu )-g_{ij}^{1} (\lambda,\mu)\right|^{2} d\mu  \le \delta_{ij}(\lambda), i,j=\overline{1,N}\right\},\]
where spectral densities $V( \lambda,\mu ),U( \lambda,\mu ),G_{1} ( \lambda,\mu )$ are known and fixed. The class $ D_V^U$ describes the "strip'' model of stochastic fields, $D_{2\delta}$ describes "$\delta$-district" in the space $L_2$ of the given spectral density  $G_1(\lambda,\mu)$.

From the condition $0\in \partial \Delta _{D} (F^{0} ,G^{0} )$ we find the following equations which determine the least favourable spectral densities for these given sets of admissible spectral densities
For the first pair ${D_{V}^{U}} ^{1}\times D_{2\delta}^{1}$ we have equations
\[((A_{T} (\lambda,\mu ))^{\top} G^{0} (\lambda,\mu )+(C_{T}^{0} (\lambda,\mu))^{\top} )^{*}((A_{T} (\lambda,\mu ))^{\top} G^{0} (\lambda,\mu )+(C_{T}^{0}(\lambda,\mu ))^{\top} )=\]
\begin{equation} \label{eq_5_1}
=(\alpha ^{2}(\lambda) +\gamma _{1}(\lambda,\mu )+\gamma _{2} (\lambda,\mu ))(F^{0} (\lambda,\mu )+G^{0} (\lambda,\mu))^{2},
\end{equation}
\[((A_{T} (\lambda,\mu ))^{\top} F^{0} (\lambda,\mu )-(C_{T}^{0} (\lambda,\mu))^{\top} )^{*}((A_{T} (\lambda,\mu ))^{\top} F^{0} (\lambda,\mu )-(C_{T}^{0}(\lambda,\mu ))^{\top})=\]
\begin{equation} \label{eq_5_2}
=\beta ^{2}(\lambda) {\mathrm{Tr}}\, (G^{0}(\lambda,\mu )-G_{1} (\lambda,\mu ))(F^{0} (\lambda,\mu )+G^{0} (\lambda,\mu ))^{2},
\end{equation}
\begin{equation} \label{eq_5_3}
\frac{1}{2\pi } \int _{-\infty }^{\infty }\left|{\mathrm{Tr}}\, (G^0(\lambda,\mu )-G_{1} (\lambda,\mu ))\right|^{2} d\mu  =\delta(\lambda) ,
\end{equation}
where $\gamma _{1} (\lambda,\mu )\le 0$ and $\gamma _{1} (\lambda,\mu )=0$ if ${\mathrm{Tr}}\,F^{0} (\lambda,\mu )> {\mathrm{Tr}}\,  V(\lambda,\mu )$, $\gamma _{2} (\lambda,\mu )\ge 0$ and $\gamma _{2} (\lambda,\mu )=0$ if $ {\mathrm{Tr}}\,F^{0}(\lambda,\mu )< {\mathrm{Tr}}\,  U(\lambda,\mu).$

For the second pair ${D_{V}^{U}} ^{2}\times D_{2\delta}^{2}$ we have equations
\[((A_{T} (\lambda,\mu ))^{\top} G^{0} (\lambda,\mu )+(C_{T}^{0} (\lambda,\mu))^{\top} )^{*}((A_{T} (\lambda,\mu ))^{\top} G^{0} (\lambda,\mu )+(C_{T}^{0}(\lambda,\mu ))^{\top} )=\]
\[=(F^{0} (\lambda,\mu)+G^{0} (\lambda,\mu ))\left\{(\alpha_{k}^{2}(\lambda) +\gamma _{1k} (\lambda,\mu )+\gamma _{2k}(\lambda,\mu ))\delta _{kl} \right\}_{k,l=1}^{N}\times\]
\begin{equation}\label{eq_5_4}
\times (F^{0} (\lambda,\mu )+G^{0}(\lambda,\mu )),
\end{equation}
\[((A_{T} (\lambda,\mu ))^{\top} F^{0} (\lambda,\mu )-(C_{T}^{0} (\lambda,\mu))^{\top} )^{*}((A_{T} (\lambda,\mu ))^{\top} F^{0} (\lambda,\mu )-(C_{T}^{0}(\lambda,\mu ))^{\top})=\]
\[=(F^{0} (\lambda,\mu )+G^{0}(\lambda,\mu ))\left\{\beta _{k}^{2}(\lambda) (g_{kk}^{0} (\lambda,\mu )-g_{kk}^{1}(\lambda,\mu ))\delta _{kl} \right\}_{k,l=1}^{N}\times\]
\begin{equation} \label{eq_5_5}
\times (F^{0} (\lambda,\mu )+G^{0}(\lambda,\mu )),
\end{equation}
\begin{equation} \label{eq_5_6}
\frac{1}{2\pi } \int _{-\infty }^{\infty }\left|g_{kk}^0 (\lambda,\mu )-g_{kk}^{1} (\lambda,\mu )\right|^{2} d\mu  =\delta _{k}(\lambda),\; k=\overline{1,N},
\end{equation}
where $\gamma _{1k} (\lambda,\mu )\le 0$ and $\gamma _{1k} (\lambda,\mu )=0$ if $f_{kk}^{0} (\lambda,\mu )>v_{kk} (\lambda,\mu ),$ $\gamma _{2k} (\lambda,\mu )\ge 0$ and $\gamma _{2k} (\lambda,\mu )=0$ if $f_{kk}^{0} (\lambda,\mu )<u_{kk} (\lambda,\mu)$.

For the third pair ${D_{V}^{U}} ^{3}\times D_{2\delta}^{3}$ we have equations
\[((A_{T} (\lambda,\mu ))^{\top} G^{0} (\lambda,\mu )+(C_{T}^{0} (\lambda,\mu))^{\top} )^{*}((A_{T} (\lambda,\mu ))^{\top} G^{0} (\lambda,\mu )+(C_{T}^{0}(\lambda,\mu ))^{\top} )=\]
\begin{equation}\label{eq_5_7}
=(\alpha ^{2}(\lambda) +\gamma'_{1} (\lambda,\mu )+\gamma'_{2} (\lambda,\mu))(F^{0} (\lambda,\mu )+G^{0} (\lambda,\mu ))(B_{1})^{\top}(F^{0} (\lambda,\mu)+G^{0} (\lambda,\mu )),
\end{equation}
\[((A_{T} (\lambda,\mu ))^{\top} F^{0} (\lambda,\mu )-(C_{T}^{0} (\lambda,\mu))^{\top} )^{*}((A_{T} (\lambda,\mu ))^{\top} F^{0} (\lambda,\mu )-(C_{T}^{0}(\lambda,\mu ))^{\top})=\]
\begin{equation} \label{eq_5_8}
=\beta ^{2}(\lambda) \left\langle B_{2},G^{0} (\lambda,\mu )-G_{1} (\lambda,\mu )\right\rangle(F^{0} (\lambda,\mu )+G^{0} (\lambda,\mu ))^2,
\end{equation}
\begin{equation} \label{eq_5_9}
\frac{1}{2\pi } \int _{-\infty }^{\infty }\left|\left\langle B_{2} ,G^0(\lambda,\mu )-G_{1} (\lambda,\mu )\right\rangle \right|^2d\mu  =\delta(\lambda),
\end{equation}
where $\gamma'_{1}( \lambda,\mu )\le 0$ and $\gamma'_{1} ( \lambda,\mu )=0$ if $\langle B_{1},F^{0} ( \lambda,\mu \rangle > \langle B_{1},V( \lambda,\mu ) \rangle,$ $\gamma'_{2}( \lambda,\mu )\ge 0$ and $\gamma'_{2} ( \lambda,\mu )=0$ if $\langle B_{1} ,F^{0} ( \lambda,\mu \rangle < \langle B_{1} ,U( \lambda,\mu ) \rangle$.

For the fourth pair ${D_{V}^{U}} ^{4}\times D_{2\delta}^{4}$ we have equations
\[((A_{T} (\lambda,\mu ))^{\top} G^{0} (\lambda,\mu )+(C_{T}^{0} (\lambda,\mu))^{\top} )^{*}((A_{T} (\lambda,\mu ))^{\top} G^{0} (\lambda,\mu )+(C_{T}^{0}(\lambda,\mu ))^{\top} )=\]
\begin{equation}\label{eq_5_10}
=(F^{0} (\lambda,\mu )+G^{0} (\lambda,\mu ))(\vec{\alpha }(\lambda)\cdot \vec{\alpha}(\lambda)^{*}+\Gamma _{1} (\lambda,\mu )+\Gamma _{2} (\lambda,\mu ))(F^{0} (\lambda,\mu)+G^{0} (\lambda,\mu ))
\end{equation}
\[((A_{T} (\lambda,\mu ))^{\top} F^{0} (\lambda,\mu )-(C_{T}^{0} (\lambda,\mu))^{\top} )^{*}((A_{T} (\lambda,\mu ))^{\top} F^{0} (\lambda,\mu )-(C_{T}^{0}(\lambda,\mu ))^{\top})=\]
\[=(F^{0} (\lambda,\mu )+G^{0} (\lambda,\mu ))\left\{\beta _{ij}(\lambda) (g_{ij}^{0} (\lambda,\mu )-g_{ij}^{1}(\lambda,\mu ))\right\}_{i,j=1}^{N}\times\]
\begin{equation} \label{eq_5_11}
\times (F^{0} (\lambda,\mu )+G^{0} (\lambda,\mu )),
\end{equation}
\begin{equation} \label{eq_5_12}
\frac{1}{2\pi } \int _{-\infty }^{\infty }\left|g_{ij}^0 (\lambda,\mu )-g_{ij}^{1} (\lambda,\mu )\right|^{2} d\mu  =\delta_{ij}(\lambda),\; i,j=\overline{1,N}.
\end{equation}
where $\Gamma _{1} (\lambda,\mu )\le 0$ and $\Gamma _{1} (\lambda,\mu )=0$ if $F^{0}(\lambda,\mu )>V(\lambda,\mu )$, $\Gamma _{2} (\lambda,\mu )\ge 0$ and $\Gamma _{2} (\lambda,\mu )=0$ if $F^{0}(\lambda,\mu )<U(\lambda,\mu )$.

The following theorem and corollaries hold true.
\begin{theorem}
Let the minimality condition hold true. The least favorable spectral densities  $F^0(\lambda,\mu)$, $G^0(\lambda,\mu)$  in the classes $D=D_{V}^U \times D_{2\delta}$ for the optimal linear estimation of the functional $A_K\vec{\xi}$ where $K=R\times[0,T]$ are determined by relations
(\ref{eq_5_1}) -- (\ref{eq_5_3}) for the first pair ${D_{V}^{U}} ^{1}\times D_{2\delta}^{1}$ of sets of admissible spectral densities;
(\ref{eq_5_4}) -- (\ref{eq_5_6}) for the second pair ${D_{V}^{U}} ^{2}\times D_{2\delta}^{2}$ of sets of admissible spectral densities;
(\ref{eq_5_7}) -- (\ref{eq_5_9}) for the third pair ${D_{V}^{U}} ^{3}\times D_{2\delta}^{3}$ of sets of admissible spectral densities;
(\ref{eq_5_10}) -- (\ref{eq_5_12}) for the fourth pair ${D_{V}^{U}} ^{4}\times D_{2\delta}^{4}$ of sets of admissible spectral densities;
constrained optimization problem (\ref{extrem}) and restrictions  on densities from the corresponding classes $D=D_{V}^U \times D_{2\delta}$.  The minimax-robust spectral characteristic of the optimal estimate of the functional $A_K\vec{\xi}$ is determined by the formula (\ref{sphar}).
\end{theorem}

\begin{corollary}
Let the minimality condition hold true. The least favorable spectral densities $F^{0}(\lambda,\mu)$ in the classes ${D_{V}^{U}} ^{k}$, $k=1,2,3,4$, for the optimal linear estimation of the functional  $A_K\vec{\xi}$,  which depends on the unknown values of the field  $\vec{\xi}(s,t)$ based on observations of the field $\vec{\xi}(s,t)$ at points $(s,t)\in R^2\setminus(R\times[0,T])$, are determined by the following  equations, respectively,
\begin{equation}
((C_{T}^{0}(\lambda,\mu) )^{\top} )^{*}\cdot(C_{T}^{0}(\lambda,\mu) )^{\top}=(\alpha^{2}(\lambda) +\gamma _{1} (\lambda,\mu )+\gamma _{2} (\lambda,\mu )) (F^{0} (\lambda,\mu ))^{2},
\end{equation}
\[((C_{T}^{0}(\lambda,\mu) )^{\top} )^{*}\cdot(C_{T}^{0}(\lambda,\mu) )^{\top}=\]
\begin{equation}
=F^{0} (\lambda,\mu )\left\{(\alpha_{k}^{2}(\lambda) +\gamma _{1k} (\lambda,\mu )+\gamma _{2k} (\lambda,\mu ))\delta _{kl} \right\}_{k,l=1}^{N}F^{0} (\lambda,\mu ),
\end{equation}
\[((C_{T}^{0}(\lambda,\mu) )^{\top} )^{*}\cdot(C_{T}^{0}(\lambda,\mu) )^{\top}=\]
\begin{equation}
=(\alpha^{2}(\lambda) +\gamma'_{1}(\lambda,\mu )+\gamma'_{2}(\lambda,\mu )) F^{0} (\lambda,\mu )(B_1)^\top F^{0} (\lambda,\mu ),
\end{equation}
\[((C_{T}^{0}(\lambda,\mu) )^{\top} )^{*}\cdot(C_{T}^{0}(\lambda,\mu) )^{\top}=\]
\begin{equation}
=F^{0} (\lambda,\mu )(\vec{\alpha}(\lambda)\cdot \vec{\alpha}(\lambda)^{*}+\Gamma _{1} (\lambda,\mu )+\Gamma _{2} (\lambda,\mu ))F^{0} (\lambda,\mu ),
\end{equation}
constrained optimization problem (\ref{extrem2}) and restrictions  on densities from the class $D_{V}^{U}$. The minimax spectral characteristic of the optimal estimate of the functional $A_K\vec{\xi}$ is determined by the formula (\ref{sphar1}).
\end{corollary}

\begin{corollary}
Let the minimality condition hold true. The least favorable spectral densities $F^{0}(\lambda,\mu)$ in the classes $D_{2\delta}^{k}$, $k=1,2,3,4$, for the optimal linear estimation of the functional  $A_K\vec{\xi}$,  which depends on the unknown values of the field  $\vec{\xi}(s,t)$ based on observations of the field $\vec{\xi}(s,t)$ at points $(s,t)\in R^2\setminus(R\times[0,T])$, are determined by the following  equations, respectively,
\begin{equation}
((C_{T}^{0}(\lambda,\mu) )^{\top} )^{*}\cdot(C_{T}^{0}(\lambda,\mu) )^{\top}=\beta ^{2}(\lambda) {\mathrm{Tr}}\, (F^{0} (\lambda,\mu )-G_{1} (\lambda,\mu ))(F^{0} (\lambda,\mu ))^{2} ,
\end{equation}
\[((C_{T}^{0}(\lambda,\mu) )^{\top} )^{*}\cdot(C_{T}^{0}(\lambda,\mu) )^{\top}=\]
\begin{equation}
=F^{0} (\lambda,\mu )\left\{\beta _{k}(\lambda)^{2} (f_{kk}^{0} (\lambda,\mu)-g_{kk}^{1} (\lambda,\mu ))\delta _{kl} \right\}_{k,l=1}^{N} F^{0} (\lambda,\mu ),
\end{equation}
\begin{equation}
((C_{T}^{0}(\lambda,\mu) )^{\top} )^{*}\cdot(C_{T}^{0}(\lambda,\mu) )^{\top}=\beta ^{2}(\lambda) \left\langle B_{2} ,F^{0} (\lambda,\mu )-G_{1} (\lambda,\mu)\right\rangle (F^{0} (\lambda,\mu ))^2,
\end{equation}
\[((C_{T}^{0}(\lambda,\mu) )^{\top} )^{*}\cdot(C_{T}^{0}(\lambda,\mu) )^{\top}=\]
\begin{equation}
=F^{0} (\lambda,\mu )\left\{\beta _{ij}(\lambda) (f_{ij}^{0} (\lambda,\mu)-g_{ij}^{1} (\lambda,\mu ))\right\}_{i,j=1}^{N} F^{0} (\lambda,\mu ),
\end{equation}
constrained optimization problem (\ref{extrem2}) and the following restrictions  on densities from the corresponding classes  $D_{2\delta}^{k}$, $k=1,2,3,4$, respectively,
\begin{equation}
\frac{1}{2\pi } \int _{-\infty }^{\infty }\left|{\mathrm{Tr}}\, (F^0(\lambda,\mu )-G_{1} (\lambda,\mu ))\right|^{2} d\mu  =\delta(\lambda) ,
\end{equation}
\begin{equation}
\frac{1}{2\pi } \int _{-\infty }^{\infty }\left|f_{kk}^0 (\lambda,\mu )-g_{kk}^{1} (\lambda,\mu )\right|^{2} d\mu  =\delta _{k}(\lambda),\; k=\overline{1,N},
\end{equation}
\begin{equation}
\frac{1}{2\pi } \int _{-\infty }^{\infty }\left|\left\langle B_{2} ,F^0(\lambda,\mu )-G_{1} (\lambda,\mu )\right\rangle \right|^2d\mu  =\delta(\lambda),
\end{equation}
\begin{equation}
\frac{1}{2\pi } \int _{-\infty }^{\infty }\left|f_{ij}^0 (\lambda,\mu )-g_{ij}^{1} (\lambda,\mu )\right|^{2} d\mu  =\delta_{ij}(\lambda),\; i,j=\overline{1,N}.
\end{equation}
The minimax spectral characteristic of the optimal estimate of the functional $A_K\vec{\xi}$ is determined by the formula (\ref{sphar1}).
\end{corollary}

\section{Conclusions}

In the article we propose methods of the mean-square optimal linear estimation of the functionals which depend on the unknown values of the multidimensional homogeneous random field based on observed data of the field with noise.  
Under condition of spectral uncertainty, which means that the spectral densities of the fields are not exactly known while a set of admissible spectral densities is given, the minimax method is applied.
Formulas for determination the least favorable spectral densities and the minimax-robust spectral characteristics of the optimal estimates of the functionals are proposed for some specific classes of admissible spectral densities. Analogous results are derived for the case of observations of the field without noise.

\end{document}